\documentclass{amsart}

\usepackage[leqno]{amsmath}
\usepackage{latexsym,amsthm,amsfonts,amssymb,xypic}
\usepackage{hyperref}

\newcommand{\numberseries}{\bfseries}   

\newlength{\thmtopspace}                
\newlength{\thmbotspace}                
\newlength{\thmheadspace}               
\newlength{\thmindent}                  

\newtheoremstyle{bfupright head,upright body}
                {\thmtopspace}{\thmbotspace}
                {\upshape}{\thmindent}{\bfseries}{.}{\thmheadspace}
                {{\numberseries \thmnumber{#2\;}}\thmnote{#3}}

\newtheoremstyle{fixed bf head,slanted body}
                {\thmtopspace}{\thmbotspace}{\slshape}
                {\thmindent}{\bfseries}{.}{\thmheadspace}
                {{\numberseries \thmnumber{#2\;}}\thmname{#1}\thmnote{ (#3)}}

\newtheoremstyle{fixed bf head,upright body}
                {\thmtopspace}{\thmbotspace}{\upshape}
                {\thmindent}{\bfseries}{.}{\thmheadspace}
                {{\numberseries \thmnumber{#2\;}}\thmname{#1}\thmnote{ (#3)}}

\theoremstyle{fixed bf head,slanted body}
\newtheorem{res}{}[section]          
\newtheorem{thm}[res]{Theorem}          \newtheorem*{thm*}{Theorem}

\newtheorem{lem}[res]{Lemma}         

\theoremstyle{bfupright head,upright body}
\newtheorem{bfhpg}[res]{}               \newtheorem*{bfhpg*}{}

\newlength{\thmlistleft}        
\newlength{\thmlistright}       
\newlength{\thmlistpartopsep}   
\newlength{\thmlisttopsep}      
\newlength{\thmlistparsep}      
\newlength{\thmlistitemsep}     

\newcounter{prt}
\newenvironment{prt}{\begin{list}{\upshape (\alph{prt})}%
    {\usecounter{prt}%
      \setlength{\leftmargin}{\thmlistleft}%
      \setlength{\labelwidth}{\thmlistleft}%
      \setlength{\rightmargin}{\thmlistright}%
      \setlength{\partopsep}{\thmlistpartopsep}%
      \setlength{\topsep}{\thmlisttopsep}%
      \setlength{\parsep}{\thmlistparsep}%
      \setlength{\itemsep}{\thmlistitemsep}}}%
  {\end{list}}%

\newenvironment{prf*}[1][Proof]{%
  \begin{proof}[\bf #1]
    \setcounter{equation}{0}
    }
  {\end{proof}
}

\newcommand{\pgref}[1]{\ref{#1}}
\newcommand{\lemref}[2][Lemma~]{#1\pgref{lem:#2}}
\renewcommand{\eqref}[1]{(\pgref{eq:#1})}
\newcommand{\thmcite}[2][?]{\cite[Thm.~#1]{#2}}
\newcommand{\corcite}[2][?]{\cite[Cor.~#1]{#2}}
\newcommand{\prpcite}[2][?]{\cite[Prop.~#1]{#2}}
\newcommand{\lemcite}[2][?]{\cite[Lem.~#1]{#2}}

\def\urltilda{\kern -.15em\lower .7ex\hbox{\~{}}\kern .04em} 

\newcommand{\setof}[3][\mspace{1mu}]{\{#1#2 \mid #3#1\}}
\newcommand{\ZZ}{\mathbb{Z}}
\newcommand{\qtext}[1]{\quad\text{#1}\quad}
\newcommand{\qqtext}[1]{\qquad\text{#1}\qquad}
\newcommand{\qand}{\qtext{and}}
\newcommand{\qqand}{\qqtext{and}}
\newcommand{\deq}{\:=\:}
\newcommand{\dge}{\:\ge\:}
\newcommand{\dle}{\:\le\:}
\newcommand{\dis}{\:\is\:}
\newcommand{\dqis}{\:\qis\:}
\newcommand{\gra}{\alpha}
\newcommand{\grb}{\beta}
\newcommand{\grf}{\varphi}
\newcommand{\mfm}{\mathfrak{m}}
\newcommand{\mfp}{\mathfrak{p}}
\newcommand{\is}{\cong}
\newcommand{\qis}{\simeq}
\newcommand{\lra}{\longrightarrow}
\newcommand{\xra}[2][]{\xrightarrow[#1]{\;#2\;}}
\newcommand{\Rhat}{\widehat{R}}
\newcommand{\Shat}{\widehat{S}}
\newcommand{\mapdef}[4][\rightarrow]{\nobreak{#2\colon #3 #1 #4}}
\newcommand{\dmapdef}[4][\lra]{\nobreak{#2\colon #3\:#1\:#4}}
\newcommand{\Cone}[1]{\nobreak{\operatorname{Cone}#1}}
\renewcommand{\H}[2][]{\operatorname{H}_{#1}(#2)}
\newcommand{\Shift}[2][]{\mathsf{\Sigma}^{#1}{#2}}
\newcommand{\SpecR}{\operatorname{Spec}R}
\newcommand{\dptR}{\operatorname{depth}R}
\newcommand{\dimR}{\operatorname{dim}R}
\newcommand{\Spec}[1]{\operatorname{Spec}#1}
\newcommand{\E}[2][R]{\operatorname{E}_{#1}(#2)}
\newcommand{\wdt}[2][R]{\operatorname{width}_{#1}#2}
\newcommand{\dpt}[2][R]{\operatorname{depth}_{#1}#2}
\newcommand{\id}[2][R]{\operatorname{id}_{#1}#2}
\newcommand{\pd}[2][R]{\operatorname{pd}_{#1}#2}
\newcommand{\Gfd}[2][R]{\operatorname{Gfd}_{#1}#2}
\newcommand{\Gid}[2][R]{\operatorname{Gid}_{#1}#2}
\newcommand{\Hom}[3][R]{\operatorname{Hom}_{#1}(#2,#3)}
\newcommand{\RHom}[3][R]{\operatorname{\mathbf{R}Hom}_{#1}(#2,#3)}
\newcommand{\tp}[3][R]{\nobreak{#2\otimes_{#1}#3}}
\newcommand{\tpp}[3][R]{(\tp[#1]{#2}{#3})}
\newcommand{\Ltp}[3][R]{\nobreak{#2\otimes_{#1}^{\mathbf{L}}#3}}
\newcommand{\Ltpp}[3][R]{(\Ltp[#1]{#2}{#3})}
\newcommand{\Cat}[2]{{\mathsf{#2}}(#1)}
\newcommand{\D}[1][R]{\Cat{#1}{D}}
 \newcommand{\A}[1][R]{\Cat{#1}{A}}
 \newcommand{\B}[1][R]{\Cat{#1}{B}}

\newcommand{\mfn}{\mathfrak{n}}
\newcommand{\mfq}{\mathfrak{q}}
\newcommand{\grfhat}{\hat{\grf}}

\numberwithin{equation}{res}

\begin{document}

\title[A bass equality for modules finite over homomorphisms]{A Bass
  equality for Gorenstein injective dimension of modules finite over
  homomorphisms}

\author[L.\,W. Christensen]{Lars Winther Christensen}

\address{Department of Mathematics and Statistics, Texas Tech
  University, Lubbock, TX 79409, U.S.A.}

\email{lars.w.christensen@ttu.edu}

\urladdr{http://www.math.ttu.edu/\urltilda lchriste}

\author[Dejun Wu]{Dejun Wu}

\address{Department of Applied Mathematics, Lanzhou University of
  Technology, Lanzhou 730050, China}

\email{wudj@lut.cn}

\thanks{L.W.C.\ was partly supported by Simons Foundation collaboration
  grant 428308; D.W.\ was partly supported by NSF of China grants
  11761047 and 11861043. The paper was written during D.W.'s year-long
  visit to Texas Tech University; the hospitality of the TTU
  Department of Mathematics and Statistics is acknowledged with
  gratitude.}

\date{30 April 2019}



\begin{abstract}
  Let $R\to S$ be a local ring homomorphism and $N$ a finitely
  generated $S$-module. We prove that if the Gorenstein injective
  dimension of $N$ over $R$ is finite, then it equals the depth of
  $R$.
\end{abstract}

\maketitle

\thispagestyle{empty}


\section*{Introduction}
\label{sec:Introduction}

\noindent
The homological theory of modules over commutative noetherian rings
comes out particularly elegant for finitely generated modules. One
way to relax this finiteness condition---without sacrificing
elegance---is to settle for finite generation over some noetherian,
but otherwise arbitrary, extension ring. This theme has been
systematically explored for at least fifteen years. As part of that
effort, this short paper answers an open question in Gorenstein
homological algebra.

In this paper a ring means a commutative noetherian ring.  Let $R$ and
$S$ be local rings with unique maximal ideals $\mfm$ and $\mfn$,
respectively. A ring homomorphism,
\begin{equation*}
  \dmapdef{\grf}{R}{S}\:,
\end{equation*}
is called \emph{local} if $\grf(\mfm) \subseteq \mfn$ holds. Given
such a homomorphism, every $S$-module is an $R$-module via $\grf$; a
finitely generated $S$-module is, when considered as an $R$-module,
said to be \emph{finite over $\grf$.} It is evident from the condition
$\grf(\mfm) \subseteq \mfn$ that $\mfm N \ne N$ holds for every module
$N\ne 0$ that is finite over $\grf$. This is an extension of
Nakayama's lemma for finitely generated modules, and the theme that
modules finite over $\grf$ behave much like finitely generated
$R$-modules was systematically explored by Avramov, Iyengar, and
Miller~\cite{AIM-06}. The first theorem in their study
 is the \emph{Bass Equality,} $\id{N} = \dptR$,
which holds if $N$ is finite over $\grf$ and of finite injective
dimension over $R$.  We extend this result with
\begin{thm*}
  Let $\mapdef{\grf}{R}{S}$ be a local ring homomorphism and $N\ne 0$
  a module finite over~$\grf$. If $N$ has finite Gorenstein injective
  dimension over $R$, then one has
  \begin{equation*}
    \Gid{N} \deq \dptR\:.
  \end{equation*}
\end{thm*}

\noindent
The statement here is a special case of Theorem~\ref{0}; it provides a
positive answer to Question 6.2 in the survey \cite{CFH-11} by
Christensen, Foxby, and Holm. Motivation for this question comes,
beyond the Bass Equality \thmcite[2.1]{AIM-06} cited above, from the
similar equality for finitely generated $R$-modules of finite
Gorenstein injective dimension, see Khatami, Tousi, and Yassemi
\corcite[2.5]{KTY-09}, and from the the \emph{Auslander--Bridger
  Equality,} $\Gfd{N} = \dptR - \dpt[R]{N}$, which by work of
Christensen and Iyengar~\cite{LWCSIn07} holds if $N$ is finite over
$\grf$ and of finite Gorenstein flat dimension over $R$.


\section{Preliminaries}

\noindent
The proof of the main result uses derived functors
on derived categories. Our notation is standard, and to not overload
this short paper we refer the reader to the appendix in \cite{lnm} for
unexplained notation.

Let $R$ be a ring, by an $R$-\emph{complex} we mean a complex of
$R$-modules. The derived category over $R$ is denoted $\D$. We say
that a complex $X$ has \emph{bounded homology} if $\H[i]{X}=0$
holds for $|i|\gg 0$. To capture the homological extent of a complex,
set
\begin{equation*}
  \inf{X} \deq \inf{\setof{i\in\ZZ}{\H[i]{X} \ne 0}} \qand 
  \sup{X} \deq \sup{\setof{i\in\ZZ}{\H[i]{X} \ne 0}}\:.
\end{equation*}
We write $\Gid{X}$ and $\Gfd{X}$ for the Gorenstein injective
dimension and Gorenstein flat dimension of an $R$-complex.  For a
complex with $\H{X} = 0$ it is standard to set $\inf{X} = \infty$,
$\sup{X} = -\infty$, and $\Gid{X} = -\infty = \Gfd{X}$.

We recall the main results from a paper by Christensen, Frankild, and
Holm~\cite{CFH-06}.

\begin{bfhpg}[The Bass category]
  \label{Bass}
  Let $R$ be a ring with a dualizing complex $D$. An $R$-complex $X$
  with bounded homology has finite Gorenstein injective dimension if
  and only if it belongs to the \emph{Bass category} $\B$; that is,
  if and only if the complex $\RHom{D}{X}$ has bounded homology,
  and the canonical morphism
  \begin{equation*}
    \dmapdef{\grb_D^X}{\Ltp{D}{\RHom{D}{X}}}{X}
  \end{equation*}
  is an isomorphism in $\D$.
\end{bfhpg}

\begin{bfhpg}[The Auslander category]
  \label{Auslander}
  Let $R$ be a ring with a dualizing complex $D$. An $R$-complex $X$
  with bounded homology has finite Gorenstein flat dimension if and
  only if it belongs to the \emph{Auslander category} $\A$; that
  is, if and only if the complex $\Ltp{D}{X}$ has bounded homology,
  and the canonical morphism
  \begin{equation*}
    \dmapdef{\gra_D^X}{X}{\RHom{D}{\Ltp{D}{X}}}
  \end{equation*}
  is an isomorphism in $\D$.
\end{bfhpg}

The next two lemmas slightly improve standard results \lemcite[(3.2.9)]{lnm}.

\begin{lem}
  \label{lem:A}
  Let $Q \to R$ be a ring homomorphism. Assume that $R$ has a
  dualizing complex and let $I$ be an injective $Q$-module. An
  $R$-complex $X$ belongs to $\A$ only if $\Hom[Q]{X}{I}$ belongs
  to $\B$, and the converse holds if $I$ is faithfully injective.
\end{lem}

\begin{prf*}
  Let $D$ be a dualizing complex for $R$.  By adjointness there is an
  isomorphism
  \begin{equation}
    \label{eq:1}
    \tag{$*$}
    \Hom[Q]{\Ltp{D}{X}}{I} \qis \RHom{D}{\Hom[Q]{X}{I}}
  \end{equation}
  in $\D$. It accounts for the horizontal isomorphism in the
  commutative diagram
  \begin{equation}
    \label{eq:2}
    \tag{$\Box$}
    \begin{gathered}
      \xymatrix{ \Hom[Q]{\RHom{D}{\Ltp{D}{X}}}{I}
        \ar[r]^-{\Hom[]{\gra_D^X}{I}}
        & \Hom[Q]{X}{I}\\
        \Ltp{D}{\Hom[Q]{\Ltp{D}{X}}{I}} \ar[u]^-\qis
        \ar[r]_-\qis & \Ltp{D}{\RHom{D}{\Hom[Q]{X}{I}}}
        \ar[u]^-{\grb_D^{\Hom[]{X}{I}}} }
    \end{gathered}
  \end{equation}
  and the vertical isomorphism is Hom evaluation; see Christensen and
  Holm~\prpcite[2.2(ii)]{LWCHHl09}.  If $X$ belongs to $\A$, then
  $\Hom[Q]{X}{I}$ has bounded homology by injectivity of $I$, the
  complex $\RHom{D}{\Hom[Q]{X}{I}}$ has bounded homology by \eqref{1},
  and $\grb_D^{\Hom[]{X}{I}}$ is an isomorphism by \eqref{2}; that is,
  $\Hom[Q]{X}{I}$ belongs to $\B$. Conversely, if $I$ is faithfully
  injective and $\Hom[Q]{X}{I}$ belongs to $\B$, then $X$ has bounded
  homology, it follows from \eqref{1} that the complex $\Ltp{D}{X}$ is
  has bounded homology and from \eqref{2} that $\gra_D^X$ is an
  isomorphism; that is, $X$ belongs to $\A$.
\end{prf*}

\begin{lem}
  \label{lem:B}
  Let $Q \to R$ be a ring homomorphism. Assume that $R$ has a
  dualizing complex and let $I$ be an injective $Q$-module. An
  $R$-complex $X$ belongs to $\B$ only if $\Hom[Q]{X}{I}$ belongs
  to $\A$, and the converse holds if $I$ is faithfully injective.
\end{lem}

\begin{prf*}
  Similar to the proof of \lemref{A}.
\end{prf*}

Another key result on Auslander and Bass categories comes from the
paper of Avramov and Foxby~\cite{LLAHBF97} in which the categories
were introduced.

\begin{bfhpg}[Regular homomorphisms]
  \label{regular}
  Let $R$ be local with maximal ideal $\mfm$ and $R \to R'$ be a flat
  local homomorphism such that the closed fiber $R'/\mfm R'$ is
  regular; such a homomorphism is called \emph{regular}. If $R$ has a
  dualizing complex, then $R'$ has a dualizing complex, see
  \cite[(2.11)]{LLAHBF97}, and by \corcite[(7.9)]{LLAHBF97} the
  next assertions hold:
  \begin{prt}
  \item An $R'$-complex belongs to $\A[R']$ if and only if it belongs to
    $\A$.
  \item An $R'$-complex belongs to $\B[R']$ if and only if it belongs to
    $\B$.
  \end{prt}
\end{bfhpg}


\section{The main result}

\noindent
We start by proving the main result in a special case, and then we
reduce the general case to the special.
Let $R$ be a local ring with maximal ideal $\mfm$.  A local ring
homomorphism $\mapdef{\grf}{R}{S}$ is said to have a \emph{regular
  factorization} if there is a commutative diagram of local ring
homomorphisms
  \begin{equation*}
    \xymatrix@=1.5pc{
      & R' \ar@{->>}[dr]\\
      R \ar[ru]^-{\dot{\grf}} \ar[rr]^{\grf} & & S
    }
  \end{equation*}
  where $\dot{\grf}$ is flat and the closed fiber $R'/\mfm R'$ is
  regular.

\begin{lem}
  \label{lem:1}
  Let $\mapdef{\grf}{R}{S}$ be a local ring homomorphism and $N$ an
  $S$-complex with bounded and degreewise finitely generated
  homology. Assume that $R$ has a dualizing complex and $\grf$ has a
  regular factorization. If $N$ has finite Gorenstein injective
  dimension over $R$, then one has
  \begin{equation*}
    \Gid{N} \deq \dptR - \inf{N}\:.
  \end{equation*}
\end{lem}

\begin{prf*}
  Let $D$ be a dualizing complex for $R$ and $R \to R' \to S$ a
  regular factorization of $\grf$. If $\H{N}=0$  the claim is
  trivial under the conventions from Section~1, so we may assume that
  $\H{N}$ is nonzero.  By \pgref{Bass} and \pgref{regular} the complex $N$
  belongs to $\B[R']$, so $g:=\Gid[R']{N}$ is finite. By
  \thmcite[6.3]{CFH-06} one has
  \begin{equation*}
    g \deq \dpt[]{R'} - \inf{N}\:.
  \end{equation*}
  Let $D'$ be a dualizing complex for $R'$, cf.~\pgref{regular}, and
  assume without loss of generality that it is normalized in the sense
  of \cite{LLAHBF97}. For every $R'$-complex $X$ with bounded and
  degreewise finitely generated homology one then has
  \begin{equation*}
    \tag{$\dagger$}
    \label{eq:dpt}
    \dpt[R']{X} \deq \inf{\RHom[R']{X}{D'}}
  \end{equation*}
  and
  \begin{equation*}
    \tag{$\dagger\dagger$}
    \label{eq:ddag}
    X \dqis \RHom[R']{\RHom[R']{X}{D'}}{D'} \qtext{in} \D[R']\:;
  \end{equation*}
  see \lemcite[(1.5.3), (2.6), and (2.7)]{LLAHBF97}.  Moreover,
  $\Gfd[R']{\RHom[R']{N}{D'}} = g$ holds by \corcite[6.4]{CFH-06}.
  Set $n = -\dimR - \inf{N}$; by \thmcite[3.1]{LWCSIn07} there is a
  distinguished triangle in $\D[R']$ of complexes with bounded and
  degreewise finitely generated homology,
  \begin{equation*}
    \tag{$\bigtriangleup$}
    \label{eq:tri}
    \RHom[R']{N}{D'} \lra P \lra H \lra \Shift{\RHom[R']{N}{D'}}\:,
  \end{equation*}
   where
  \begin{equation*}
    \label{eq:H}
    \tag{$*$}
    \pd[R']{P} \deq g \qqand \sup{H} \dle \Gfd[R']{H} \dle n\:.
  \end{equation*}
  By the  Auslander--Buchsbaum formula and \thmcite[4.1]{LWCSIn07} one has
  \begin{equation*}
    \dptR' - \dpt[R']{P} \deq g \deq \dptR' - \dpt[R']{\RHom[R']{N}{D'}}\:.
\end{equation*}
 Combined with \eqref{dpt}  and \eqref{ddag} these equalities  yield
  \begin{equation*}
    \label{eq:di}
    \tag{$**$}
    \dpt[R']{P} \deq \inf{N}\:.
  \end{equation*}

  Applying the functor $\RHom[R']{-}{D'}$ to \eqref{tri} one gets via
  \eqref{ddag} the triangle
  \begin{equation*}
    \label{eq:tritri}
    \tag{$\bigtriangledown$}
    \Shift[-1]{N} \lra \RHom[R']{H}{D'} \lra \RHom[R']{P}{D'} \lra N\:.
  \end{equation*}
  One has $\id[R']{\RHom[R']{P}{D'}} = g$; see
  \corcite[6.4]{CFH-06}. As $\dot{\grf}$ is flat, the complex
  $\RHom[R']{P}{D'}$ has finite injective dimension over $R$. By
  \corcite[8.2.2]{AIM-06} and \thmcite[4.4]{DWu15} one has
  $\id{\RHom[R']{P}{D'}} = \dptR - \inf{\RHom[R']{P}{D'}}$, which by
  \eqref{dpt} and \eqref{di} can be rewritten as
  \begin{equation*}
    \label{eq:id}
    \tag{$***$}
    \id{\RHom[R']{P}{D'}} \deq \dptR - \inf{N}\:.
  \end{equation*}
  The complex $\RHom[R']{H}{D'}$ has finite Gorenstein injective
  dimension over $R'$ by \corcite[6.4]{CFH-06} and hence over $R$;
  see \pgref{regular}. The first inequality in the next
  computation holds by \thmcite[3.3]{CFH-06}; the equality follows
  from \eqref{dpt}; the second inequality holds by the definition of
  depth; the final inequality follows from \eqref{H}.
  \begin{align*}
    \Gid{\RHom[R']{H}{D'}}  
    & \dle \dimR - \inf{\RHom[R']{H}{D'}}\\
    & \deq \dimR - \dpt[R']{H}\\
    & \dle \dimR + \sup{H}\\
    & \dle -\inf{N}\,.
  \end{align*}
  For every injective $R$-module $I$ and every $i \le \inf{N}$ one
  gets from \eqref{tritri} an exact sequence in homology
  \begin{equation*}
    \H[i]{\RHom{I}{\RHom[R']{P}{D'}}} \lra \H[i]{\RHom{I}{N}} \lra 0\:.
  \end{equation*}
  Thus, per \thmcite[3.3]{CFH-06} and  \eqref{id} one has
  \begin{equation*}
    \Gid{N} \le \id{\RHom[R']{P}{D'}} = \dptR-\inf{N}\:.
\end{equation*}
The opposite inequality $\Gid{N} \ge \dptR-\inf{N}$ holds by
\thmcite[6.3]{CFH-06}.
\end{prf*}

As is standard, we denote by $\Rhat$ and $\Shat$ the completions of
$R$ and $S$ in the topologies induced by their maximal ideals.
  The homomorphism $\mapdef{\grf}{R}{S}$ extends to a homomorphism of
  complete local rings; that is, there is a commutative diagram of
  local ring homomorphisms
  \begin{equation*}
    \xymatrix@=1.5pc{
      \Rhat \ar[r]^-\grfhat & \Shat \\
      R \ar[u] \ar[r]^-\grf & S \ar[u]
    }
  \end{equation*}
In particular, every $\Shat$-complex is an $\Rhat$-complex.

In the special case $R=S$, and $\grf$ the identity, the next result was
proved by Christensen, Frankild, and Iyengar; see Foxby and
Frankild~\thmcite[3.6]{HBFAJF07}.

\begin{lem}
  \label{lem:2}
  Let $\mapdef{\grf}{R}{S}$ be a local ring homomorphism and $N$ an
  $S$-complex with bounded and degreewise finitely generated
  homology. If $N$ has finite Gorenstein injective dimension over $R$,
  then $\tp[S]{N}{\Shat}$ has finite Gorenstein injective dimension
  over $\Rhat$.
\end{lem}

\begin{prf*}
  Let $K^S$ be the Koszul complex on a minimal set of generators for
  $\mfn$, the maximal ideal of $S$. Since the $S$-complex $\H{K}$ has
  degreewise finite length, one has $\tp[S]{\Shat}{K^S} \qis K^S$ in
  $\D[S]$. Under the flat map $S \to \Shat$ the minimal generators of
  $\mfn$ extend to a minimal set of generators for the maximal ideal
  $\hat{\mfn}$ of $\Shat$, so $\tp[S]{\Shat}{K^S}$ is the Koszul
  complex $K^{\Shat}$ on a minimal set of generators for
  $\hat{\mfn}$. Thus one has
  \begin{equation}
    \label{eq:K}
    \tag{$\diamond$}
    K^S \qis K^{\Shat}
  \end{equation}
  in $\D[S]$, and we simply denote this complex $K$.
  
  The first step is to notice that $\tp[S]{N}{K}$ has finite
  Gorenstein injective dimension over $R$.  For every element
  $x \in \mfn$ there is an exact sequence of $S$-complexes,
  \begin{equation*}
    0 \lra N \lra \Cone{x^N} \lra \Shift{N} \lra 0\:,
  \end{equation*}
  where $x^N$ is the homothety. Since $N$ and $\Shift{N}$ have finite
  Gorenstein injective dimension over $R$, so has $\Cone{x^N}$; this
  folklore fact is dual to a result of Veliche \thmcite[3.9]{OVl06} for
  Gorenstein projective dimension. Now, $\Cone{x^N}$ is isomorphic to
  $\tp[S]{N}{K(x)}$, where $K(x)$ denotes the elementary Koszul
  complex on $x$. Since $K$ is a tensor product of such elementary
  Koszul complexes, it follows that $\tp[S]{N}{K}$ has finite
  Gorenstein injective dimension over $R$.

  Set $M = \tp[S]{N}{K}$; it is an $\Shat$-complex via $K$ and, therefore, an
  $\Rhat$-complex. The second step is to prove that $M$ belongs to
  $\B[\Rhat]$. The composite $R \xra{\grf} S \lra \Shat$, called the
  semi-completion of $\grf$, has a regular factorization; see Avramov,
  Foxby, and Herzog \thmcite[(1.1)]{AFH-94}.  Let $E'$ denote the injective hull
  of the residue field of $R'$ and $E$ denote the injective hull
  $\E{k} \is \E[\Rhat]{k}$.  As $\H{M}$ has degreewise finite length
  over $\Shat$ and, therefore, over $R'$ one has
  \begin{equation*}
    M \dqis \Hom[R']{\Hom[R']{M}{E'}}{E'}\:.
  \end{equation*}
  As $\Gid{M}$ is finite, it follows from \thmcite[1.7]{LWCSSW10} that
  the $R'$-complex
  \begin{align*}
    \RHom{R'}{M} &\dqis \RHom{R'}{\Hom[R']{\Hom[R']{M}{E'}}{E'}}\\
                 &\dqis \Hom[R']{\tp{R'}{\Hom[R']{M}{E'}}}{E'}
  \end{align*}
  has finite Gorenstein injective dimension. The $R'$-complex
  $\tp{R'}{\Hom[R']{M}{E'}}$ then has finite Gorenstein flat dimension
  by \lemref{A}. As $R'$ is faithfully flat over $R$, it follows from
  \thmcite[1.8]{LWCSSW10} that the complex $\Hom[R']{M}{E'}$ has finite
  Gorenstein flat dimension over $R$. By another application of
  the same result the complex $\tp{\Rhat}{\Hom[R']{M}{E'}}$ has finite
  Gorenstein flat dimension over $\Rhat$, whence it belongs to
  $\A[\Rhat]$. By \lemref{A} the dual complex
  \begin{align*}
    \Hom[\Rhat]{\tp{\Rhat}{\Hom[R']{M}{E'}}}{E} \dis \Hom{\Hom[R']{M}{E'}}{E}
  \end{align*}
  belongs to $\B[\Rhat]$. As $E$ is faithfully injective, it follows
  from \lemref{A} that the complex $\Hom[R']{M}{E'}$ belongs to $\A[\Rhat]$ and
  hence to $\A[R']$; see \pgref{regular}. By \lemref{B} the complex
  $M$ now belongs to $\B[R']$ and hence to $\B[\Rhat]$.

  To finish the proof we now prove that $\tp[S]{N}{\Shat}$ belongs to
  $\B[\Rhat]$, cf.~\ref{Bass}.  First notice that by \eqref{K} and
  associativity of the tensor product one has
  \begin{equation*}
    \label{eq:4}
    \tag{$\diamond\diamond$}
    \tp[S]{N}{K} \dis
    \tp[S]{N}{\tpp[\Shat]{\Shat}{K}} \dis \tp[\Shat]{\tpp[S]{N}{\Shat}}{K}\:.
  \end{equation*}
  By \eqref{4} and an application of
  tensor evaluation \prpcite[2.2(v)]{LWCHHl09} one gets
  \begin{equation*}
    \tag{$\ddag\ddag$}
    \begin{aligned}
      \RHom[\Rhat]{D}{\tp[S]{N}{K}}
      &\dqis \RHom[\Rhat]{D}{\tp[\Shat]{\tpp[S]{N}{\Shat}}{K}}\\
      &\dqis \tp[\Shat]{\RHom[\Rhat]{D}{\tp[S]{N}{\Shat}}}{K}\:.
    \end{aligned}
  \end{equation*}
  As $\tp[S]{N}{K}$ belongs to $\B[\Rhat]$, the complex
  $\tp[\Shat]{\RHom[\Rhat]{D}{\tp[S]{N}{\Shat}}}{K}$ has bounded
  homology, so $\RHom[\Rhat]{D}{\tp[S]{N}{\Shat}}$ has bounded
  homology. This follows from work of Foxby and
  Iyengar~\cite[1.3]{HBFSIn03}; indeed, as the $\Rhat$-complex $D$ and
  the $\Shat$-complex $\tp[S]{N}{K}$ have degreewise finitely
  generated homology, it follows from \lemcite[1.3.2]{AIM-06} that the $\Shat$-complex
  $\RHom[\Rhat]{D}{\tp[S]{N}{\Shat}}$ has degreewise finitely
  generated homology. There is a
  commutative diagram in $\D[\Rhat]$,
  \begin{equation*}
    \xymatrix@C=4pc{
      \Ltp[\Rhat]{D}{\RHom[\Rhat]{D}{\tp[S]{N}{K}}}
      \ar[d]_-\qis \ar[r]^-{\grb_D^{\tp[]{N}{K}}}_-\qis 
      & \tp[S]{N}{K} \ar[d]^-\qis \\ 
      \tp[\Shat]{\Ltpp[\Rhat]{D}{\RHom[\Rhat]{D}{\tp[S]{N}{\Shat}}}}{K}
      \ar[r]^-{\tp[]{\grb_D^{\tp[]{N}{\Shat}}}{K}} 
      & \tp[\Shat]{\tpp[S]{N}{\Shat}}{K}
    }
  \end{equation*}
  where the right-hand vertical isomorphism is \eqref{4}, and the
  left-hand vertical isomorphism follows by tensor evaluation
  \prpcite[2.2(v)]{LWCHHl09} and associativity of the tensor
  product. It follows that
  ${\tp[\Shat]{\grb_D^{\tp[S]{N}{\Shat}}}{K}}$ is an isomorphism; that
  is, the mapping cone
  \begin{equation*}
    \Cone{{\tpp[\Shat]{\grb_D^{\tp[S]{N}{\Shat}}}{K}}} 
    \dqis \tp[\Shat]{(\Cone{\grb_D^{\tp[S]{N}{\Shat}})}}{K}
  \end{equation*}
  is acyclic. As $\grb_D^{\tp[S]{N}{\Shat}}$ per
  \lemcite[1.3.2]{AIM-06} is a morphism of $\Shat$-complexes with
  degreewise finitely generated homology, it follows from
  \cite[1.3]{HBFSIn03} that the complex $\Cone{\grb_D^{\tp[S]{N}{\Shat}}}$ is
  acyclic, whence $\grb_D^{\tp[S]{N}{\Shat}}$ is an isomorphism in $\D[\Rhat]$, and
  $\tp[S]{N}{\Shat}$ belongs to $\B[\Rhat]$. 
\end{prf*}

The main result, which we can now prove, compares to \thmcite[4.1 and
Cor.~4.8]{LWCSIn07}.

\begin{thm}
  \label{0}
  Let $\mapdef{\grf}{R}{S}$ be a local ring homomorphism and $N$ an
  $S$-complex with bounded and degreewise finitely generated
  homology. If $N$ has finite Gorenstein injective dimension over $R$,
  then one has
  \begin{align*}
    \Gid{N} & \deq \dptR - \inf{N}\\
            & \deq -\inf{\RHom{\E{k}}{N}}\\
    & \deq \Gid[\Rhat]{\tpp[S]{N}{\Shat}}\:.
  \end{align*}
\end{thm}

\begin{prf*}
  The homomorphism $\Rhat \xra{\grfhat} \Shat$ has a regular
  factorization; see \thmcite[(1.1)]{AFH-94}. By \lemref{2} the
  $\Rhat$-complex $\tp[S]{N}{\Shat}$ has finite Gorenstein injective
  dimension, and it has bounded and degreewise finite homology over
  $\Shat$, so \lemref{1} yields
  \begin{equation*}
    \Gid[\Rhat]{\tpp[S]{N}{\Shat}} \deq \dpt[]{\Rhat} - \inf{\tpp[S]{N}{\Shat}}\:.
  \end{equation*}
  There are equalities $\dpt[]{\Rhat} = \dptR$ and
  $\inf{\tpp[S]{N}{\Shat}} = \inf{N}$; the latter holds by faithful
  flatness of $\Shat$ over $S$. Moreover, one has
  \begin{equation*}
    \Gid{N} \dge \dptR - \inf{N} \deq -\inf{\RHom{\E{k}}{N}}
  \end{equation*}
  by \thmcite[6.3]{CFH-06} and \corcite[6.5]{LWCHHl09}, so it is
  sufficient to prove that the inequality
  $\Gid{N} \le \Gid[\Rhat]{\tpp[S]{N}{\Shat}}$ holds.  By
  \thmcite[2.2]{LWCSSW10} one has
  \begin{equation*}
    \label{eq:ss}
    \tag{\S}
    \begin{aligned}
      \Gid{N} &\deq \sup{\setof{\dpt[]{R_\mfp} -
          \wdt[R_\mfp]{N_\mfp}}{\mfp\in\SpecR}}\quad\text{and}\\
      \Gid[\Rhat]{\tpp[S]{N}{\Shat}} &\deq
      \sup{\setof{\dpt[]{\Rhat_\mfq} -
          \wdt[\Rhat_\mfq]{\tpp[S]{N}{\Shat}_\mfq}}{\mfq\in\Spec{\Rhat}}}\:.
    \end{aligned}
  \end{equation*}
  For $\mfp\in \SpecR$ choose $\mfq\in \Spec{\Rhat}$ minimal over
  $\mfp\Rhat$, the local homomorphism $R_\mfp \to \Rhat_\mfq$ is flat
  with artinian closed fiber, whence one has
  $\dptR_\mfp = \dpt[]{\Rhat_\mfq}$; see e.g.\
  \prpcite[(2.8)]{AFH-94}. In the next computation the first and
  fourth equalities hold by the definition of width, the second holds
  by faithful flatness of $\Shat$ over $S$, and the last
  holds as $R_\mfp \to \Rhat_\mfq$ is a
  local homomorphism; see Wu and Kong~\lemcite[3.6]{DWuFKn18}.
  \begin{align*}
    \wdt[R_\mfp]{N_\mfp} 
    & \deq \inf{\tpp{(R_\mfp/\mfp R_\mfp)}{N}}\\
    & \deq \inf{\tp[S]{\tpp{(R_\mfp/\mfp R_\mfp)}{N}}{\Shat}}\\
    & \deq \inf{\tp{(R_\mfp/\mfp R_\mfp)}{\tpp[S]{N}{\Shat}}}\\
    & \deq \wdt[R_\mfp]{\tpp[S]{N}{\Shat}_\mfp}\\
    & \deq \wdt[\Rhat_\mfq]{\tpp[S]{N}{\Shat}_\mfq}\:.
  \end{align*}
  For every prime ideal $\mfp$ in $R$ there is thus a prime ideal
  $\mfq$ in $\Rhat$ with
  \begin{equation*}
    \dpt[]{R_\mfp} - \wdt[R_\mfp]{N_\mfp} = \dpt[]{\Rhat_\mfq} -
    \wdt[\Rhat_\mfq]{\tpp[S]{N}{\Shat}_\mfq}
  \end{equation*}
  so the desired inequality is immediate from \eqref{ss}.
\end{prf*}


\section*{Acknowledgment}
\noindent 
We thank the anonymous referee for suggestions and comments that helped
us improve the exposition.



\def\soft#1{\leavevmode\setbox0=\hbox{h}\dimen7=\ht0\advance \dimen7
  by-1ex\relax\if t#1\relax\rlap{\raise.6\dimen7
  \hbox{\kern.3ex\char'47}}#1\relax\else\if T#1\relax
  \rlap{\raise.5\dimen7\hbox{\kern1.3ex\char'47}}#1\relax \else\if
  d#1\relax\rlap{\raise.5\dimen7\hbox{\kern.9ex \char'47}}#1\relax\else\if
  D#1\relax\rlap{\raise.5\dimen7 \hbox{\kern1.4ex\char'47}}#1\relax\else\if
  l#1\relax \rlap{\raise.5\dimen7\hbox{\kern.4ex\char'47}}#1\relax \else\if
  L#1\relax\rlap{\raise.5\dimen7\hbox{\kern.7ex
  \char'47}}#1\relax\else\message{accent \string\soft \space #1 not
  defined!}#1\relax\fi\fi\fi\fi\fi\fi}
  \providecommand{\MR}[1]{\mbox{\href{http://www.ams.org/mathscinet-getitem?mr=#1}{#1}}}
  \renewcommand{\MR}[1]{\mbox{\href{http://www.ams.org/mathscinet-getitem?mr=#1}{#1}}}
  \providecommand{\arxiv}[2][AC]{\mbox{\href{http://arxiv.org/abs/#2}{\sf
  arXiv:#2 [math.#1]}}} \def\cprime{$'$}
\providecommand{\bysame}{\leavevmode\hbox to3em{\hrulefill}\thinspace}
\providecommand{\MR}{\relax\ifhmode\unskip\space\fi MR }
\providecommand{\MRhref}[2]{%
  \href{http://www.ams.org/mathscinet-getitem?mr=#1}{#2}
}
\providecommand{\href}[2]{#2}

\end{document}